\documentclass[10pt]{amsart}
\usepackage{a4}
\usepackage{amssymb}
\usepackage{amscd}
\newtheorem{thm}{Theorem}

\newtheorem{cor}[thm]{Corollary}
\newtheorem{lem}[thm]{Lemma}
\newtheorem{prop}[thm]{Proposition}
\newtheorem{defn}[thm]{Definition}

\newtheorem{rem}[thm]{Remark}
\newtheorem*{tthm}{Theorem}

\numberwithin{thm}{section}
\numberwithin{equation}{section}

\newcommand{\norm}[1]{\left\Vert#1\right\Vert}
\newcommand{\abs}[1]{\left\vert#1\right\vert}

\newcommand{\F}{\mathcal{F}}
\newcommand{\G}{\mathcal{G}}
\newcommand{\Hi}{\mathcal{H}}
\newcommand{\N}{\mathbb{N}}

\newcommand{\U}{\mathcal{U}}

\newcommand{\Om}{\Omega}
\newcommand{\om}{\omega}
\newcommand{\Le}{\mathcal{L}}
\newcommand{\R}{\mathcal{R}}

\begin{document}
\title{Type and cotype of operator spaces}
\author{Hun Hee Lee}
\address{Department of Mathematical Sciences, Seoul National University
          San56-1 Shinrim-dong Kwanak-gu Seoul 151-747, Republic of Korea}
\email{lee.hunhee@gmail.com, bbking@amath.kaist.ac.kr}

\begin{abstract}
We consider two operator space versions of type and cotype, namely $S_p$-type, $S_q$-cotype and type $(p,H)$, cotype $(q,H)$ 
for a homogeneous Hilbertian operator space $H$ and $1\leq p \leq 2 \leq q\leq \infty$, generalizing ``$OH$-cotype 2" of G. Pisier.
We compute type and cotype of some Hilbertian operator spaces and $L_p$ spaces, and we investigate the relationship between 
a homogeneous Hilbertian space $H$ and operator spaces with cotype $(2,H)$. As applications we consider
operator space versions of generalized little Grothendieck's theorem and Maurey's extension theorem in terms of these new notions.
\end{abstract}
\maketitle

\section{Introduction}
Type and cotype plays an important role in Banach space theory. Thus, it is natural to expect operator space analogue of type and cotype.
Actually there has been several attempts to define type and cotype in the operator space category.
In \cite{P4} G. Pisier defined $OH$-cotype 2 and M. Junge (chapter 4. of \cite{J}) studied a variant of this notion, namely cotype $(2, R+C)$,
where $OH$ and $R+C$ imply the operator Hilbert space and operator space sum of the row and column Hilbert spaces, respectively.
In this paper we are going to give two different definitions of type and cotype of operator spaces, namely 
$S_p$-type, $S_q$-cotype and type $(p,H)$, cotype $(q,H)$ for a homogeneous Hilbertian operator space $H$, which 
are both generalizations of Pisier's ``$OH$-cotype 2". 

In order to get a satisfactory theory we need to focus on two aspects.  The first one is about how big the cotype 2 class is. 
Note that cotype 2 class in Banach space category includes all (noncommutative) $L_p$ ($1\leq p \leq 2$) spaces. 
The second one is about the possibility of applications such as generalized little Grothendieck's theorem, 
Maurey's extension theorem and Kwapie\'{n}'s theorem. Although these new notions still do not promise satisfactory results in both ways, 
each definitions have their own pros. For example, type $(p,H)$ and cotype $(q,H)$ of $L_p$ spaces behave well for some good choice of $H$, 
and both notions allow corresponding applications mentioned above.
Besides, cotype $(2,H)$ gives us an insight about the relationship between a homogeneous Hilbertian space $H$ 
and an operator space $E$ which has cotype $(2,H)$. More precisely, it is known that (\cite{J, P4}) 
$S_1$ (the trace class on $\ell_2$) have cotype $(2, R+C)$ but not cotype $(2,OH)$. 
Thus, it is natural to be interested which $H$ is best among all $H$, which $S_1$ has cotype $(2,H)$. 
This question will be answered later for all $L_p$ ($1\leq p \leq 2$) spaces and the meaning of ``best" will be clarified.

Now let us discuss our approach more precisely. Recall that a Banach space $X$ is called (gaussian) type $p$ ($1\leq p \leq 2$) 
if and only if there exist a constant $C>0$ such that
\begin{eqnarray} \label{def-type-p}
\pi_{p',2}(v^*) \leq C\cdot \ell^*(v)
\end{eqnarray}
for $\frac{1}{p} + \frac{1}{p'} =1$ and every $v: X \rightarrow \ell^n_2$, $n\in \mathbb{N}$ 
and (gaussian) cotype $q$ ($2\leq q \leq \infty$) if and only if there exist a constant $C'>0$ such that
\begin{eqnarray}\label{def-cotype-q}
\pi_{q,2}(u) \leq C'\cdot \ell(u)
\end{eqnarray}
for every $u: \ell^n_2 \rightarrow X$ and $n\in \mathbb{N}$.
Here, $\pi_{q,2}(\cdot)$ is the $(q,2)$-summing norm defined by
$$\pi_{q,2}(T:X\rightarrow Y) = \sup\Big\{ \frac{(\sum_k \norm{Tx_k}^q )^{\frac{1}{q}}}
{\norm{\sum_k x_k\otimes e_k}_{X \otimes_{\lambda}\ell_2}} \Big\},$$ 
where $\otimes_{\lambda}$ is the injective tensor product of Banach spaces,
and $\ell(u)$ and $\ell^*(v)$ are $\ell$-norm and
adjoint $\ell$-norm, respectively, defined by
$$\ell(u) := \Big[\int_{\Omega} \norm{\sum^n_{k=1}g_k(\omega)ue_k}^2_X
dP(\omega)\Big]^{\frac{1}{2}}$$ for an i.i.d. gaussian variables
$\{g_k\}$ on a probability space $(\Omega, P)$ and 
$$\ell^*(v) := \sup\{\text{\rm tr}(vu) | u:\ell^n_2 \rightarrow X, \ell(u)\leq 1 \}.$$

Pisier's definition of $OH$-cotype 2 is as follows.
{\it An operator space $E$ is said to have ``$OH$-cotype 2" if there is a constant $C>0$ such that
for all $v : E \rightarrow \ell^n_2$, we have
\begin{eqnarray}\label{def-OH-cotype-2}
\ell^*(v) \leq C\pi_{2,oh}(v),
\end{eqnarray}
} 
where $\pi_{2,oh}(v)$ is the $(2,oh)$-summing norm defined by 
$$\pi_{2,oh}(v) = \sup\Big\{ \frac{(\sum_k \norm{vx_k}^2 )^{\frac{1}{2}}}{\norm{\sum_k x_k\otimes e_k}_{E\otimes_{\text{min}}OH}} \Big\}$$
and $\otimes_{\text{min}}$ is the injective tensor product of operator spaces.
Note that this definition is based on the trace dual formulation of (\ref{def-cotype-q}). 
Thus, in order to extend these notions for the general exponents $q \geq 2$ we need to consider trace dual version of (\ref{def-OH-cotype-2}). 
However, unlike in the Banach space case we have the problem that $\pi_{2,oh}$ is not self-dual. We can resolve this difficulty by observing that (Proposition 6.2 in \cite{P3}) $$\pi_{2,oh}(v : E\rightarrow \ell^n_2) = \pi^o_2(v : E\rightarrow OH_n),$$ 
where $\pi^o_2(\cdot)$ is the completely 2-summing norm defined by 
$$\pi^o_2(T : E\rightarrow F) = \sup\Big\{ \frac{\norm{(Tx_{ij})}_{S_2(F)}}{\norm{\sum_{i,j} x_{ij}\otimes e_{ij}}_{E\otimes_{\text{min}}S_2}} \Big\}$$
and $S_r(F)$ ($1\leq r \leq \infty$) implies vector valued Schatten classes introduced in \cite{P3}.
Then since $\pi^o_2$ is self-dual (which will be checked later) we can reformulate (1.5) as follows.
$E$ is ``OH-cotype $2$" if there is a constant $C>0$ such that
for all $u : OH_n \rightarrow E$, we have
\begin{eqnarray}
\pi^o_2(u) \leq C \cdot \ell(u).
\end{eqnarray}
Now it is easy to extend cotype notion to $q\geq 2$ case by replacing $\pi^o_2(u)$ into 
the completely $(q,2)$-summing norm $\pi^o_{q,2}(u)$ defined by 
$$\pi^o_{q,2}(T : E\rightarrow F) := \sup\Big\{ \frac{\norm{(Tx_{ij})}_{S_q(F)}}{\norm{\sum_{i,j} x_{ij}\otimes e_{ij}}_{E\otimes_{\text{min}}S_2}} \Big\}$$
and we will call it as $S_q$-cotype. $S_p$-type can be defined similarly.

There is another approach using approximation numbers. This can be done in a more general context. 
Let $H$ be a homogeneous Hilbertian operator space, i.e. $H$ is isometric to a Hilbert space and for every $u : H \rightarrow H$ 
we have $\norm{u}_{cb}  = \norm{u}$. Then we can define $\pi_{2,H}(v)$ replacing $OH$ into $H$ and use it in the definition of cotype $(2,H)$.
In order to assure that $\pi_{2,H}(\cdot)$ is actually a norm we need to assume that $H$ is ``subquadratic"
i.e. for all orthogonal projections $\{P_i\}^n_{i=1}$ in $H$ with $I_H = P_1 + \cdots + P_n$ we have 
$$\norm{x}^2_{B(\ell_2)\otimes_{\min}H} \leq \sum^n_{i=1}\norm{I_{B(\ell_2)}\otimes P_i (x)}^2_{B(\ell_2)\otimes_{\min}H}$$
for any $x\in B(\ell_2)\otimes H$. (See p.82 of \cite{P2}) 

$E$ is called ``cotype $(2,H)$" if there is a constant $C>0$ such that for all $u : \ell^n_2 \rightarrow E$, we have
\begin{eqnarray}
\pi^*_{2,H}(u) \leq C \cdot \ell(u),
\end{eqnarray}
where $\pi^*_{2,H}$ implies the trace dual of $\pi_{2,H}$.
Now we recall the equivalence between $\pi_{q,2}(u)$ and $(\sum_k a_k(u)^q)^{\frac{1}{q}}$ for 
$u : \ell_2 \rightarrow X$, (Corollary 19.7 of \cite{TJ}) where $a_k(\cdot)$ is the $k$-th approximation number 
defined by $$a_k(u) = \inf \{ \norm{u-v} : v \in B(\ell_2,X), \text{rk}(v) < k\}.$$
Since we do not have appropriate $(q,2)$-extension of $\pi^*_{2,H}$ we use $\ell_q$-sum of c.b. approximation numbers of the map 
$u :H^*_n \rightarrow E$, where $H^*_n$ be the $n$-dimensional version of $H^*$. See section \ref{sec-type(p,H)-cotype(q,H)} for the details. 
We will call it as cotype $(q,H)$, and type $(p,H)$ can be defined similarly. Note that $S_2$-type and $S_2$-cotype is equivalent to type $(2,OH)$ and cotype $(2,OH)$, respectively.

The behavior of $S_q$-cotype of $L_p$ spaces are quite different from that of Banach space case.
However, the behavior of cotype $(2,H)$ depends on $H$. More precisely we have the following.
\begin{tthm}
Let $1\leq p \leq 2$ and $\mu$ be a $\sigma$-finite measure.
\begin{itemize}
\item[(1)]
$S_p$ is cotype $(2,H)$ if and only if the formal identity $$id : RC[p] \rightarrow H$$ is completely bounded.
\item[(2)]
$L_p(\mu)$ is cotype $(2,H)$ if and only if the formal identity $$id : RC[p'] \rightarrow H$$ is completely bounded.
\end{itemize}
\end{tthm}
Note that $RC[r] = [R\cap C, R+C]_{\frac{1}{r}}$ for $1\leq r \leq \infty$, where $\cap$ and $[\cdot , \cdot]_{\frac{1}{r}}$ 
imply the intersection and the complex interpolation in the category of operator spaces, respectively. 
Thus we can say that $RC[p]$ (resp. $RC[p']$) is the best (in the above sense) homogeneous Hilbertian operator space, 
which $S_p$ (resp. $L_p(\mu)$) has cotype $(2,H)$. 

With these notions of type and cotype we can consider several applications.
The first one is an operator space analogue of ``generalized little Grothendieck's theorem".
(See \cite{DPR, M1} for the Banach space case)

\begin{tthm}
Every bounded linear map from $C(K)$ into $S_q$-cotype space is completely $(q,2)$-summing for a compact set $K$ and
$2\leq q<\infty$.
\end{tthm}

The second one is an operator space analogue of ``Maurey's extension theorem". (See \cite{M2} for the Banach space case)
We say that an operator space $H$ is perfectly Hilbertian if $H$ is a homogeneous Hilbertian operator space and $H$ and $H^*$ are subquadratic.
(See section 8. of \cite{P2}) Note that $R[p] = [R,C]_{\frac{1}{p}}$ and $C[p] = [C,R]_{\frac{1}{p}}$
are perfectly Hilbertian, and $RC[p]$ are completely isomorphic to perfectly Hilbertian operator spaces. Then we have the following.
\begin{tthm}
Let $E$ and $F$ be operator spaces with type $(2,H)$ and cotype $(2,H^*)$,
respectively, for a perfectly Hilbertian operator space $H$. Then there is a constant $C>0$ 
such that for any subspace $G \subseteq E$ and any bounded linear map $u : G \rightarrow F$ 
we have an extension $$\tilde{u} : E \rightarrow F \,\, \text{with}\,\, \gamma_H(\tilde{u}) 
\leq C\norm{u}.$$ 
\end{tthm}
Recall that $\gamma_H(\tilde{u}) = \inf \{\norm{A}_{cb}\norm{B}_{cb}\},$ 
where the infimum runs over all possible factorization $\tilde{u}: E \stackrel{A}{\longrightarrow}
H(I) \stackrel{B}{\longrightarrow} F$ for some index set $I$. We need ``perfectness" of $H$ to ensure that $\gamma_H(\cdot)$ is actually a norm.
By the Remark in p.82 of \cite{P2} $H(I)$ is well defined for any index set $I$.
As a corollary we get operator space versions of ``Kwapie\'{n}'s theorem". See \cite{K} for classical Banach space case and 
\cite{GP2} for another operator space case.

Note that there is a different notion of type and cotype of operator spaces by J. Garcia-Cuerva and J. Parcet
using quantized orthonormal systems (\cite{GP1, GP2, Par}). 
We will see how $S_2$-type and $S_2$-cotype is related to the type 2 and cotype 2 in \cite{GP2} at the end of section 2.

This paper is organized as follows. In section \ref{sec-SptypeSqcotype}, we define $S_p$-type and $S_q$-cotype of operator spaces 
and develop their basic theory. As examples, we estimate $S_p$-type and $S_q$-cotype of $R[p] = [R,C]_{\frac{1}{p}}$, $C[p] = [C,R]_{\frac{1}{p}}$
and $L_p$ spaces. $S_p$-type and $S_q$-cotype of commutative $L_p$ spaces are the same as in the Banach space case 
while those of $S_p$ are completely different. In section \ref{sec-type(p,H)-cotype(q,H)} we define type $(p,H)$ and cotype $(q,H)$ of operator spaces
and investigate their basic properties. As examples, we compute type $(p,H)$ and cotype $(q,H)$ of $R[p]$, $C[p]$ and $L_p$ spaces.
Note that we can recover the same result as in the Banach space case for $S_p$ with a ``good" choice of $H$.
Moreover, we investigate the relationship between a homogeneous Hilbertian space $H$ and operator spaces with cotype $(2,H)$.
In the last section we present the above two applications of type, cotype notions in this paper. 

Throughout this paper, we will assume some knowledge of operator space theory (\cite{EF, P5}), completely $p$-summing
maps (\cite{P3}), absolutely $p$-summing operators (\cite{DJT, TJ}) and vector valued noncommutative $L_p$-spaces (\cite{P3}). 
For $1\leq p \leq\infty$ $S_p$ (resp. $S^n_p$) and $S_p(E)$ (resp. $S^n_p(E)$) refer to Schatten classes on $\ell_2$ (resp. $\ell^n_2$)
and their vector valued versions. (\cite{P3})
We will assume that all $L_p$ spaces (commutative or noncommutative) and their vector valued versions are endowed with
their natural operator space structure in the sense of \cite{P3}. In this paper $H$ is reserved for a homogeneous Hilbertian operator space on $\ell_2$ we will denote its $n$-dimensional version by $H_n$. In particular, $R_n[p]$, $C_n[p]$ and $RC_n[p]$ imply $n$-dimensional versions of $R[p]$, $C[p]$ and $RC[p]$ respectively. As usual, $B(E,F)$ and $CB(E,F)$ denote the set of all bounded linear maps and all completely bounded linear maps from $E$ into
$F$, respectively. We use the symbol $a\lesssim b$ if there is a $C>0$ such that $a \leq C b$ and $a\sim b$ if $a\lesssim b$ and
$b\lesssim a$. We denote the conjugate exponent of $1\leq r\leq \infty$ by $r'$, i.e. $\frac{1}{r} + \frac{1}{r'} = 1$.

\section{$S_p$-type and $S_q$-cotype of operator spaces} \label{sec-SptypeSqcotype}

\subsection{Definition and basic properties}

As an operator space version of ``absolutely $p$-summing operators" G. Pisier introduced ``completely $p$-summing maps" in
\cite{P3} as follows.\\ {\it A linear map between operator spaces $u : E\rightarrow F$ is called ``completely $p$-summing" for
$1\leq p \leq \infty $ if $$I_{S_p}\otimes u : S_p\otimes_{\min} E \rightarrow S_p(F)$$ is a bounded map.} 
We denote $\pi^o_p(u)$ for the operator norm of $I_{S_p}\otimes u$. Similarly we define an operator space version of ``absolutely $(q,2)$-summing operators".\\ {\it A linear map between operator spaces $u : E\rightarrow F$ is
called ``completely $(q,2)$-summing" for $2 \leq q \leq \infty $ if $$I_{2,q}\otimes u : S_2\otimes_{\min} E \rightarrow S_q(F)$$
is a bounded map}, where $I_{2,q}$ is the formal identity from $S_2$ into $S_q$. We denote $\pi^o_{q,2}(u)$ for the operator norm
of $I_{2,q}\otimes u$ and $\Pi^o_{q,2}(E,F)$ for the collection of all such operators from $E$ into $F$. Now we define $S_p$-type and $S_q$-cotype.

\begin{defn}
Let $E$ be an operator space.

\begin{itemize}

\item[(1)] 
$E$ is said to have $S_p$-type ($1\leq p \leq 2$) if there is a constant $C >0$ such that
$$\pi^o_{p',2}(v^*) \leq C \cdot \ell^*(v)$$ for every $n\in \mathbb{N}$ and $v: E \rightarrow OH_n$.

\item[(2)] 
$E$ is said to have $S_q$-cotype ($2\leq q \leq \infty$) if there is a constant $C' >0$ such that
$$\pi^o_{q,2}(u) \leq C' \cdot \ell(u)$$ for every $n\in \mathbb{N}$ and $u: OH_n \rightarrow E$.
\end{itemize}

\end{defn}

We can reformulate $S_p$-type and $S_q$-cotype by comparing vector-valued Schatten class norm of $E$-valued matrices and their gaussian averages.
Let $\{g_{ij}\}$ be an re-indexing of $\{ g_i \}$.

\begin{prop}\label{prop-Sptype-gaussian-average}
Let $E$ be an operator space.

\begin{itemize}

\item[(1)] For $n \in \mathbb{N}$ and $1\leq p\leq 2$ we define $T_{S_p,n}(E)$ to be the infimum of the constant $C>0$ satisfying
\begin{eqnarray}\label{Sptype-gaussian-average}
\Big[\int_{\Om}\norm{\sum^n_{i,j=1}g_{ij}(\om)x_{ij}}^{2} dP(\om)\Big]^{\frac{1}{2}}\leq C \norm{(x_{ij})}_{S^n_p(E)}.
\end{eqnarray}
$E$ has \textbf{$S_p$-type} if and only $$T_{S_p}(E) = \sup_{n\geq 1}T_{S_p,n}(E) < \infty.$$

\item[(2)] For $n \in \mathbb{N}$ and $2\leq q\leq \infty$ we define $C_{S_q,n}(E)$ to be the infimum of the constant $C'>0$ satisfying
\begin{eqnarray}\label{Sqcotype-gaussian-average}
\norm{(x_{ij})}_{S^n_q(E)}\leq C' \Big[\int_{\Om}\norm{\sum^n_{i,j=1}g_{ij}(\om)x_{ij}}^{2} dP(\om)\Big]^{\frac{1}{2}}.
\end{eqnarray}
$E$ has \textbf{$S_q$-cotype} if and only if $$C_{S_q}(E) = \sup_{n\geq 1}C_{S_q,n}(E) < \infty.$$
\end{itemize}
\end{prop}

\begin{proof}
(1) $E$ satisfies (\ref{Sptype-gaussian-average}) if and only if there is a constant
$C>0$ such that we have $$\ell(u) \leq C
\norm{(ue_{ij})}_{S^n_p(E)}$$ for all $n\in \mathbb{N}$ and $u:
S^n_2 \rightarrow E$. By trace duality this is equivalent to
$$\norm{(v^*e_{ij})}_{S^n_{p'}(E^*)} \leq C \cdot \ell^*(v)$$ for all $n \in
\mathbb{N}$ and $v : E \rightarrow S^n_2$. Indeed, by Corollary 1.8
of \cite{P2} we have
\begin{align*}
\begin{split}
\abs{\text{tr}(vu)} & = \abs{\sum_{i,j} \langle vue_{ij}, e_{ij}
\rangle} = \abs{\sum_{i,j} \langle ue_{ij}, v^*e_{ij} \rangle}
\\ & \leq \norm{(ue_{ij})}_{S^n_p(E)}\norm{(v^*e_{ij})}_{S^n_q(E^*)}.
\end{split}
\end{align*}
Now we assume that $E$ satisfies (\ref{Sptype-gaussian-average}) and consider $v : E \rightarrow OH_n$ and $(x_{ij}) \in M_m(OH_n)$, 
$m\in \mathbb{N}$. If we set $w : S^m_2 \rightarrow OH_n, e_{ij} \mapsto x_{ij}$, then we have by 12.5 in \cite{TJ} that
\begin{align*}
\begin{split}
\norm{(v^*x_{ij})}_{S^m_{p'}(E^*)} & = \norm{(v^*we_{ij})}_{S^m_{p'}(E^*)}
\leq C \cdot \ell^*(w^*v) \\ & = C\sup\{ \abs{\text{tr}(uw^*v)} :
\ell(u : S^m_2 \rightarrow E) \leq 1 \}
\\ & \leq C\sup\{ \abs{\text{tr}(\tilde{u}v)} :
\ell(\tilde{u} : OH_n \rightarrow E) \leq 1 \}\norm{w}\\ & = C \cdot
\ell^*(v)\norm{w}_{cb} = C \cdot
\ell^*(v)\norm{(x_{ij})}_{S^m_2\otimes_{\min}OH_n}.
\end{split}
\end{align*}
The converse direction is straightforward from the above observation and the fact that 
$$\norm{(e_{ij})}_{S^n_2 \otimes_{\min} S^n_2} = \norm{I_{S^n_2}}_{cb} =1.$$

(2) Suppose $E$ satisfies (\ref{Sqcotype-gaussian-average}) and let $u: OH_n \rightarrow E$ and $(x_{ij}) \in
S^m_q(OH_n)$ for $m \in \mathbb{N}$. If we set $v: S^m_2 \rightarrow OH_n, e_{ij} \mapsto x_{ij}$, 
then we have $\norm{v} = \norm{v}_{cb} = \norm{(x_{ij})}_{S^m_2\otimes_{\min}OH_n}$. Now, we have by (12.5) in \cite{TJ}
\begin{align*}
\begin{split}
\norm{(ux_{ij})}_{S^m_q(E)} & = \norm{(uve_{ij})}_{S^m_q(E)} \leq C' \norm{\sum_{1 \leq i,j \leq m}g_{ij}(\cdot)uve_{ij}}_{L_2(\Om,E)}\\
& = C'\ell(uv) \leq C'\ell(u)\norm{v} = C'\ell(u)\norm{(x_{ij})}_{S^m_2\otimes_{\min}OH_n}.
\end{split}
\end{align*}
The converse direction is straightforward as before.

\end{proof}

\begin{rem}\label{rem-def}{ \rm
\begin{itemize}

\item[(1)] If we take diagonals of (\ref{Sptype-gaussian-average}) and (\ref{Sqcotype-gaussian-average}), then it is trivial that every
$S_p$-type (resp. $S_q$-cotype) space has type $p$ (resp. cotype $q$) as a Banach space.

\item[(2)] Instead of gaussian systems we can use the Rademacher system $\{r_i\}$ defined by $r_i(t) = \text{sign}(\sin(2^i\pi t))$,
$t\in [0,1]$ and $i=1,2,\cdots$ in the definition to get the Rademacher $S_p$-type and $S_q$-cotype. 
It is easy to check that two notions are equivalent when $1< p \leq 2 \leq q < \infty$. 
Although we don't know the equivalence for the case $p=1$ and $q = \infty$ all the calculations in this paper can be transfered 
to the Rademacher setting with the same argument. 

\item[(3)] Unlike in the Banach space case, $S_1$-type and $S_{\infty}$-cotype are no more trivial, that is, 
we have examples of operator spaces without $S_1$-type and $S_{\infty}$-cotype, respectively. 
We will see examples in Lemma \ref{lem-Sptype-estimate-Sp} and Theorem \ref{thm-subhom} in detail.
Moreover, for any operator space $E$ we have $$T_{S_1, n}(E)  \lesssim n^{\frac{1}{2}}\,\,
\text{and}\,\, C_{S_{\infty},n}(E) \lesssim n^{\frac{1}{2}}.$$
Indeed, we have
\begin{align*}
\begin{split}
\norm{\sum^n_{i,j=1}r_{ij}\otimes x_{ij}}_{L_1(\Om,E)} & =
\int_{\Om} \norm{\sum^n_{i,j=1}x_{ij}g_{ij}(\om)} dP(\om)
\\ & \leq \int_{\Om} \norm{(g_{ij}(\om))}_{S^n_\infty}\norm{(x_{ij})}_{S^n_1(E)}dP(\om) \\
& \lesssim n^{\frac{1}{2}}\norm{(x_{ij})}_{S^n_1(E)}
\end{split}
\end{align*}
by Lemma 2.3 in \cite{GP1} and Proposition 45.1 in \cite{TJ}.
The estimation for $C_{S_{\infty},n}(E)$ can be obtained by the duality below. (Proposition \ref{duality})

\item[(4)] We consider the following transforms.
$$\F_\G : f \mapsto \Big( \int_{\Om} f(t)g_{ij}(\om)dP(\om) \Big)_{i,j}\,\,\,\,\text{and}\,\,\,\,
\F^{-1}_\G : (x_{ij}) \mapsto \sum_{i,j} g_{ij}(\om)x_{ij}$$ for appropriate $f : \Om \rightarrow \mathbb{C}$
and $(x_{ij}) \in M_{\infty}$. Then $E$ has $S_p$-type ($1\leq p\leq 2$) if and only if 
$$\F^{-1}_\G \otimes I_E : S_p(E) \rightarrow \G_2(E)$$ is bounded and $E$ has $S_q$-cotype ($2\leq q \leq \infty$) if
and only if $$\F_\G \otimes I_E : \G_2(E) \rightarrow S_q(E)$$ is bounded, where $\G_r(E)$ refers to the closed linear span of
$\{g_{ij}\}\otimes E$ in $L_r(\Om,E)$ for $1\leq r < \infty$. We write $\G^n_r(E)$ ($n\in \mathbb{N}$) for the closed linear span
of $\{g_{ij}\}^n_{i,j=1}\otimes E$ in $L_r(\Om,E)$.
\end{itemize}
}
\end{rem}

$S_p$-type and $S_q$-cotype have a partial duality as follows. The proof is the same as in the Banach space case, 
so that we omit it. Note that we can include the cases $S_1$-type and $S_{\infty}$-cotype without any extra effort.
See Proposition 11.10 and 13.17 in \cite{DJT}.

\begin{prop} \label{duality} Let $E$ be an operator space, $1\leq p \leq 2$ and $n\in \mathbb{N}$.
\begin{itemize}
\item[(1)] If $E$ has $S_p$-type, then $E^*$ has $S_{p'}$-cotype
with $C_{S_{p'},n}(E^*) \leq T_{S_p,n}(E)$. \item[(2)] If $E$ has
$S_{p'}$-cotype and is K-convex as a Banach space, then $E^*$ has
$S_p$-type with $T_{S_p,n}(E^*)\leq K(E)C_{S_{p'},n}(E)$,
\end{itemize}
where $K(E)$ is the K-convexity constant of $E$ defined by the operator norm of the $E$-valued gaussian projection 
from $L_2(\Om,E)$ onto $\G_2(E)$, given by $$f\mapsto \sum_{i,j}\Big(\int_{\Om} f(t)g_{ij}(\om)dP(\om)\Big) g_{ij}.$$
\end{prop}

\subsection{Relationships to other concepts}

Now we check that the $S_2$-cotype in this paper coincide with the OH-cotype 2 in \cite{P4}. 
It can be achieved by the following trace duality of $\pi^o_2$-norm. It is well-known to experts, 
but we include the proof since we could not find the reference.

\begin{lem}\label{lem2}
Let $E$ and $F$ be operator spaces and $E$ be finite dimensional.
Then for $v:F\rightarrow E$ we have
$$(\pi^o_2)^*(v) := \sup\{\abs{{\rm tr}(vu)}|\pi^o_2(u:E\rightarrow F) \leq 1\} = \pi^o_2(v).$$
\end{lem}
\begin{proof}
Let $u : E\rightarrow F$ and $v: F\rightarrow E$. By Proposition
6.1 in \cite{P3}, we have factorizations $$u : E \stackrel{V_1}{\rightarrow} OH(I) \stackrel{T_1}{\rightarrow} F\,\,
\text{and} \,\, v : F \stackrel{V_2}{\rightarrow} OH(J) \stackrel{T_2}{\rightarrow} E$$ for some index sets $I$ and $J$ with
$$\pi^o_2(V_1),\,\, \pi^o_2(V_2) \leq 1, \,\, \norm{T_1}_{cb} \leq \pi^o_2(u)\,\, \text{and}\,\, \norm{T_2}_{cb} \leq \pi^o_2(v).$$
Then, by Proposition 6.3 in \cite{P3} we have
\begin{align*}
\begin{split}
\abs{{\rm tr}(vu)} & = \abs{{\rm tr}(T_2V_2T_1V_1)} = \abs{{\rm tr}(V_1T_2V_2T_1)}\\ &\leq \norm{V_2T_1}_{HS} \norm{V_1T_2}_{HS} = \pi^o_2(V_2T_1)\pi^o_2(V_1T_2)\\ & \leq
\pi^o_2(V_2)\norm{T_1}_{cb}\pi^o_2(V_1)\norm{T_2}_{cb}\leq \pi^o_2(v)\pi^o_2(u),
\end{split}
\end{align*}
where $\norm{\cdot}_{HS}$ implies the Hilbert-Schmidt norm. 

Thus, we get $(\pi^o_2)^*(v) \leq \pi^o_2(v)$.

For the converse inequality we consider any $\epsilon > 0$ and choose
$(x_{ij}) \in S^n_2 \otimes_{\min}F$ with $$\norm{(x_{ij})}_{S^n_2 \otimes_{\min}F} = 1 \,\, \text{and}\,\, \norm{(vx_{ij})}_{S^n_2(E)}
\geq (1-\epsilon)\pi^o_2(v).$$ Then there is $(y^*_{ij}) \in S^n_2(E^*)$ such that $$\norm{(y^*_{ij})}_{S^n_2(E^*)} = 1 \,\,
\text{and}\,\, \norm{(vx_{ij})}_{S^n_2(E)} = \abs{\left\langle (vx_{ij}), (y^*_{ij}) \right\rangle}.$$
Now we set $A : E \rightarrow S^n_2, x\mapsto (y^*_{ij}x)_{i,j}$ and $B: S^n_2 \rightarrow F, e_{ij} \mapsto x_{ij}$.
Then we get $$\pi^o_2(BA) \leq \norm{B}_{cb}\pi^o_2(A) \leq \norm{(y^*_{ij})}_{S^n_2(E)}\norm{(x_{ij})}_{S^n_2\otimes_{\min}F} \leq 1$$
by Lemma 5.14 of \cite{P3}. Thus, we have
\begin{align*}
\begin{split}
(\pi^o_2)^*(v) & \geq \abs{{\rm tr}(vBA)} = \abs{{\rm tr}(AvB)} =  \sum^n_{i,j=1}\abs{\left\langle AvBe_{ij},
e_{ij} \right\rangle}\\ & =\abs{\left\langle (vx_{ij}), (y^*_{ij}) \right\rangle} \geq (1-\epsilon)\pi^o_2(v),
\end{split}
\end{align*}
which gives us the converse inequality.
\end{proof}

\begin{cor}
The $S_2$-cotype coincides with the OH-cotype 2 in \cite{P4}
\end{cor}
\begin{proof}
By Proposition 6.2 in \cite{P3}, we have $\pi^o_2(v) = \pi_{2,oh}(v)$ for any $v : E\rightarrow OH_n$. Thus, we get the
desired conclusion by Lemma \ref{lem2} and trace duality.
\end{proof}

We end this section by providing a partial relationship between $S_2$-type and $S_2$-cotype and notions in \cite{GP2}.

Let $(\Om, P)$ be a probability space and $(\Sigma, d_{\Sigma})$
be a pair of an index set $\Sigma$ and a collection of natural
numbers indexed by $\Sigma$, $d_{\Sigma} = \{ d_{\sigma} \in
\mathbb{N}: \sigma \in \Sigma \}$. The quantized gaussian system
$\G_{\Sigma}$ with parameter $(\Sigma, d_{\Sigma})$ is the
collection of random matrices $g^\sigma =
\frac{1}{\sqrt{d_{\sigma}}}(g^\sigma_{ij}): \Om \rightarrow
M_{d_\sigma}$ indexed by $\Sigma$, where $g^\sigma_{ij}$'s are
i.i.d. gaussian random variables. We consider the following
transforms.
$$\F_{\G_{\Sigma}}(f)(\sigma) = \int_\Om f(\om)g^\sigma(\om)^*dP(\om)\,\,\,\,\text{and}\,\,\,\,
\F^{-1}_{\G_{\Sigma}} (A)(\om) = \sum_{\sigma\in\Sigma}d_\sigma
\text{tr}(A^\sigma g^\sigma(\om))$$ for appropriate $f : \Om
\rightarrow \mathbb{C}$ and $A \in \prod_{\sigma\in
\Sigma}M_{d_\sigma}$.

For $1\leq p\leq 2$, $2\leq q \leq \infty$, $\frac{1}{p} +
\frac{1}{p'} = 1$ and $\frac{1}{q} + \frac{1}{q'} = 1$ we say that
an operator space $E$ has Banach $\G_{\Sigma}$-type $p$ if
$$\sup_{\text{finite}\,
\Gamma \subseteq \Sigma} \norm{\F^{-1}_{\G_{\Sigma}} \otimes
I_E}_{\Le_{p}(\Gamma, E)\rightarrow L_{p'}(\Om, E)} <\infty$$ and
that $E$ has Banach $\G_{\Sigma}$-cotype $q$ if

$$\sup_{\text{finite}\, \Gamma \subseteq
\Sigma}\norm{\F_{\G_{\Sigma}} \otimes I_E}_{L^\Gamma_{q'}(\Om, E)
\rightarrow \Le_{q}(\Gamma, E)} <\infty$$ where
$L^\Gamma_{q'}(\Om, E)$ is the closed linear span of
$\{g^\sigma_{ij} : \sigma\in \Gamma\}\otimes E$ in $L_{q'}(\Om,
E)$,
$$\Le_{r}(\Gamma, E) = \{ A\in \prod_{\sigma\in
\Gamma}M_{d_\sigma}\otimes E: \norm{A}_{\Le_{r}(\Gamma, E)} =
{\Big( \sum_{\sigma\in \Gamma}d_\sigma
\norm{A^\sigma}^r_{S^{d_\sigma}_{r}(E)}\Big)}^\frac{1}{r}<\infty\}$$
for $1 \leq r < \infty$ and
$$\Le_{\infty}(\Gamma, E) = \{ A\in
\prod_{\sigma\in \Gamma}M_{d_\sigma}\otimes E:
\norm{A}_{\Le_{\infty}(\Gamma, E)} = \sup_{\sigma\in
\Gamma}\norm{A^\sigma}_{S^{d_\sigma}_{\infty}(E)}<\infty\}.$$ For
the details and the natural operator space structure on
$\Le_{r}(\Gamma, E)$, see \cite{GP1, P3}.

\begin{prop} Let $E$ be an operator space and $\G_\Sigma$ be the quantized gaussian system with parameter
$(\Sigma, d_{\Sigma})$. Suppose that $d_{\Sigma}$ is unbounded. Then $E$ has gaussian $S_2$-type if and only if 
it has Banach $\G_\Sigma$-type 2 and $E$ has gaussian $S_2$-cotype if and only if it has Banach $\G_\Sigma$-cotype 2.
\end{prop}
\begin{proof}
Let $\Gamma$ be a finite subset of $\Sigma$ and $A (=(A^{\sigma}))\in \Pi_{\sigma \in \Gamma}M_{d_{\sigma}}\otimes E$.
If we set $$B = \oplus_{\sigma \in \Gamma}\sqrt{d_{\sigma}}A^{\sigma} \in S^n_{\infty}(E)$$ for 
$n = {\sum_{\sigma\in\Gamma}d_{\sigma}}$, then we get
$$\F^{-1}_{\G_{\Sigma}}(A)(\om) = \sum_{\sigma\in\Gamma}d_\sigma \text{tr}(A^\sigma g^\sigma(\om)) =
\sum_{\sigma\in\Gamma}\sqrt{d_\sigma} \text{tr}(A^\sigma (g^\sigma_{ij}(\om)))= \F^{-1}_{\G}(B)(\om)$$ and
$$\norm{A}_{\Le_2(\Gamma,E)} = \Big[ \sum_{\sigma\in\Gamma}d_\sigma\norm{A^\sigma}^2_{S^\sigma_2(E)} \Big]^\frac{1}{2} = \norm{B}_{S^n_2(E)}.$$ Conversely, for any $B \in S^n_{\infty}(E)$ we choose $\sigma_0 \in \Sigma$ with $d_{\sigma_0} > n$ and set 
$A (=(A^{\sigma}))\in \Pi_{\sigma \in \Sigma}M_{d_{\sigma}}\otimes E$ by $A^{\sigma_0} =d_{\sigma_0}^{-\frac{1}{2}} B \oplus 0$ 
and $A^{\sigma} = 0$ elsewhere. Then we also get $\F^{-1}_{\G_{\Sigma}}(A) = \F^{-1}_{\G}(B)$ and 
$\norm{A}_{\Le_2(\Sigma,E)} = \norm{B}_{S^n_2(E)}.$ Thus, we get the desired result.
\end{proof}

\subsection{$S_p$-type and $S_q$-cotype of $R[p]$ and $C[p]$}\label{sec-exam-RpCp}

In case of Hilbertian spaces the gaussian average of vector valued matrix is simple to calculate so that we can determine $S_p$-type
and $S_q$-cotype in some concrete cases.

\begin{thm} \label{compute1}
Let $1 \leq p \leq \infty$ and $\frac{1}{p} + \frac{1}{p'} = 1$. Then $R[p]$ (resp. $C[p])$ has $S_{\min \{p, p'\}}$-type 
and $S_{\max \{p, p'\}}$-cotype and does not have $S_r$-type nor $S_{r'}$-cotype for $\min \{p, p'\} < r \leq 2$ and 
$\frac{1}{p} + \frac{1}{p'} = 1$.
\end{thm}
\begin{proof}
Note that $R$ and $R_n$ are isometric to $OH$ and $OH_n$, respectively. Thus, we have that
\begin{align*}
\begin{split}\label{cotype_condition}
\text{$R$ has $S_q$-cotype}\,\, & \text{$\Leftrightarrow \F_\R
\otimes I_R :
\G_2(R) \rightarrow S_q(R)$ is bounded} \\
& \text{$\Leftrightarrow \F_\R \otimes id : \G_2(OH) \rightarrow
S_q(R)$ is bounded}\\
&  \text{$\Leftrightarrow I_{2,q}\otimes id : S_2(OH) \rightarrow
S_q(R)$ is bounded}
\end{split}
\end{align*}
$$\text{$\Leftrightarrow I^n_{2,q}\otimes id_n :
S^n_2(OH_n) \rightarrow S^n_q(R_n)$ is uniformly bounded for
all}\,\, n \in \mathbb{N},$$ where $id$, $id_n$, $I_{2,q}$ and
$I^n_{2,q}$ are corresponding formal identities.

First, we consider the case $q = \infty$. For $(x_{ij})\in S^n_{\infty}(R_n)$, $x_{ij} = \sum^n_{k=1}x^k_{ij}e_{1k}$ we have
$$\norm{(x_{ij})}_{S^n_2(OH_n)} = \Big(\sum^n_{i,j,k =1}\abs{x^k_{ij}}^2\Big)^{\frac{1}{2}}$$ and by considering
$(x_{ij})$ as a $n^2 \times n^2$-matrix
\begin{align*}
\begin{split}
\norm{(x_{ij})}_{S^n_{\infty}(R_n)} & = \sup\Big\{
\Big[\sum^n_{i=1}\Big(\sum^n_{j,k=1}x^k_{ij}a_{jk}\Big)^2\Big]^\frac{1}{2}\Big
| \sum^n_{j,k=1}\abs{a_{jk}}^2 = 1\Big\}.
\end{split}
\end{align*}

However, we have
\begin{align*}
\begin{split}
\Big[\sum^n_{i=1}\Big(\sum^n_{j,k=1}x^k_{ij}a_{jk}\Big)^2\Big]^\frac{1}{2}
& \leq
\Big[\sum^n_{i=1}\Big(\sum^n_{j,k=1}\abs{x^k_{ij}}^2\Big)\cdot\Big(\sum^n_{j,k=1}
\abs{a_{jk}}^2\Big)\Big]^\frac{1}{2}\\ & = \Big(\sum^n_{i,j,k
=1}\abs{x^k_{ij}}^2\Big)^{\frac{1}{2}}\Big(\sum^n_{j,k=1}
\abs{a_{jk}}^2\Big)^\frac{1}{2},
\end{split}
\end{align*}
and consequently $$\norm{(x_{ij})}_{S^n_{\infty}(R_n)} \leq \norm{(x_{ij})}_{S^n_2(OH_n)}.$$ Thus, $R$ has $S_{\infty}$-cotype
with $C_{S_{\infty}}(R) = 1$, and we can similarly show that $C$ has $S_{\infty}$-cotype with
$C_{S_{\infty}}(C) = 1$. Since $R$ and $C$ are $K$-convex as Banach spaces, $R$ and $C$ have $S_1$-type by duality (Proposition
\ref{duality}). Since $R[2]$ (resp. $C[2]$) is completely isometric to $OH$, it has $S_2$-type and $S_2$-cotype. Thus by complex
interpolation, $R[p]$ (resp. $C[p])$ has $S_{\min \{p, p'\}}$-type and $S_{\max \{p, p'\}}$-cotype.

Now suppose $2\leq q < p \leq \infty$ and consider $S_q(R[p])$. By Theorem 1.1 of \cite{P3}, 
$$S_q(R[p]) \cong C[q]\otimes_h R[p] \otimes_h R[q]$$ completely isometrically under the mapping
$$e_{ij}\otimes x \mapsto e_{i1}\otimes x\otimes e_{1j},$$ where $\otimes_h$ refers to the Haggerup tensor product.
Note that by the commutation property of the Haggerup tensor product with respect to complex interpolation we have 
the following completely isometric isomorphisms.
\begin{align*}
\begin{split}
R[p] \otimes_h R[q] & \cong [R[p]\otimes_h R, R[p]\otimes_h C]_{\frac{1}{q}}\\
& \cong [[R\otimes_h R, C\otimes_h R]_{\frac{1}{p}}, [R\otimes_h C,
C\otimes_h C]_{\frac{1}{p}}]_{\frac{1}{q}}.
\end{split}
\end{align*}
Since $C\otimes_h R$ (resp. $R\otimes_h C$) is completely isometric to $S_{\infty}$ (resp. $S_1$) and $R\otimes_h R$ and
$C\otimes_h C$ are isometric to $S_2$, respectively, we get a
subspace $$F (\cong R[p] \otimes_h R[q])\,\,\text{of}\,\, S_q(R[p]) \,\, \text{isometric to}\,\, S_r$$
under the mapping $e_{1j}\otimes e_{1i} \mapsto e_{ij}$, where $r = \frac{2pq}{pq+p-q}$. However, we have
$$(I_{2,q}\otimes id)^{-1}(e_{1j}\otimes e_{1i}) = e_{1j}\otimes e_i \in S_2(OH),$$ so that 
$$G = (I_{2,q}\otimes id)^{-1}(F) \cong S_2$$ isometrically.

Consequently, $I_{2,q}\otimes id$ cannot be bounded since $r<2$ and $(I_{2,q}\otimes id) |_G$ is nothing but the formal identity
$I_{2,r} : S_2 \rightarrow S_r$, which means $R[p]$ does not have $S_q$-cotype. We can show that $C[p] $ does not have 
$S_q$-cotype similarly, and type cases are obtained by the duality (Proposition \ref{duality}). Since $R[p'] \cong C[p]$ we get the
desired result for all $1\leq p \leq \infty$.
\end{proof}

\subsection{$S_p$-type and $S_q$-cotype of $L_p$ spaces}\label{sec-exam-Lp}

In this section we will compute $S_p$-type and $S_q$-cotype of $L_p$ spaces. Unfortunately we don't have good
behavior as in the Banach spaces cases generally. We only have the same results as in the Banach space case for $L_p$ spaces 
($1\leq p < \infty$) with respect to Type I von Neumann algebras of bounded degree. When $p = \infty$, we have very bad behavior even
in the commutative cases.

\begin{thm}\label{Thm-SemiNoncomLp}
Let $(\mathcal{M}, \mu)$ be a $\sigma$-finite measure space, $1\leq p <\infty$ and $n\in \mathbb{N}$.

\begin{itemize}

\item[(1)] $L_p(\mu, S^n_p)$ has $S_r$-type and $S_{r'}$-cotype for $r = \min\{p,2\}$ and $\frac{1}{r} + \frac{1}{r'} = 1$. 
If $L_p(\mu)$ is infinite dimensional, then it
does not have $S_r$-type nor $S_{r'}$-cotype for any $\min\{p,2\}< r\leq 2$.

\item[(2)] $L_{\infty}(\mu, S^n_{\infty})$ has $S_{\infty}$-cotype.
If $L_{\infty}(\mu)$ infinite dimensional, then it does not have
$S_1$-type nor $S_q$-cotype for any $q< \infty$.

\end{itemize}
\end{thm}
\begin{proof}
Note that $p=2$ cases are trivial. First, we consider the case $1 \leq p \leq 2$.
For $(x_{ij}) \in S_2(L_1(\mu, S^n_1))$, we have
\begin{align*}
\begin{split}
\norm{\sum_{i,j}g_{ij}\otimes x_{ij}}_{L_2(\Om,L_1(\mu, S^n_1))} & \geq
\norm{\sum_{i,j}g_{ij}\otimes x_{ij}}_{L_1(\Om,L_1(\mu, S^n_1))}\\
& = \int_{\mathcal{M}}\norm{\sum_{i,j}x_{ij}(s)\otimes g_{ij}}_{S^n_1(L_1(\Om))} d\mu(s)
\end{split}
\end{align*}
Since $\G_1$ and $S_2$ are isomorphic as Banach spaces,
$S^n_1(\G_1)$ and $S^n_1(S_2)$ are isomorphic allowing constants
depending on $n$. Indeed, we have $S^n_1(\G_1)^* = CB(\G_1,
S^n_{\infty}) \cong B(\G_1, S^n_{\infty})  \cong B(S_2,
S^n_{\infty}) \cong CB(S_2, S^n_{\infty}) = S^n_1(S_2)^*$
isomorphically. Thus, we have by Corollary 1.10 in \cite{P3} that
\begin{align*}
\begin{split}
\norm{\sum_{i,j}g_{ij}\otimes x_{ij}}_{L_2(\Om,L_1(\mu, S^n_1))}
& \gtrsim \int_{\mathcal{M}}\norm{\sum_{i,j}x_{ij}(s)\otimes e_{ij}}_{S^n_1(S_2)} d\mu(s)\\
& \geq \norm{(x_{ij})}_{S_2(L_1(\mu, S^n_1))},
\end{split}
\end{align*}
Thus, $L_1(\mu, S^n_1)$ has OH-cotype 2 and by complex interpolation with $L_2(\mu, S^n_2)$
we obtain that $L_p(\mu, S^n_p)$($1 \leq p \leq 2$) has $S_2$-cotype.

Since $S_1 \stackrel{I_{1,2}}{\longrightarrow} S_2 \stackrel{\F^{-1}_\R}{\longrightarrow} \G_2 \subseteq L_2\Om$
is a contraction, where $I_{2,\infty}$ is the corresponding formal identity, so is 
$$\F^{-1}_\R \otimes I_{L_1(\mu, S^n_1)} : S_1\otimes_{\gamma} L_1(\mu, S^n_1) \rightarrow L_2\Om\otimes_{\gamma} L_1(\mu, S^n_1),$$ 
where $\otimes_{\gamma}$ is the projective tensor product in the category of Banach spaces. 
Note that $L_2\Om\otimes_{\gamma} E \hookrightarrow L_2(\Om,E)$ contractively by the canonical embedding and 
$S_1 \otimes_{\gamma} L_1(\mu, S^n_1) \cong S_1(L_1(\mu, S^n_1))$ isomorphically. 
Indeed, we have $$(S_1 \otimes_{\gamma} L_1(\mu, S^n_1))^* = B(S_1(L_1(\mu)), S^n_{\infty}) \cong CB(S_1(L_1(\mu)), 
S^n_{\infty}) = S_1(L_1(\mu, S^n_1))^*$$ isomorphically. (allowing constants depending on $n$.)
Thus, we have a bounded map $$\F^{-1}_\R \otimes I_{L_1(\mu, S^n_1)} : S_1(L_1(\mu, S^n_1))
\rightarrow L_2(\Om, L_1(\mu, S^n_1)),$$ which implies that $L_1(\mu, S^n_1)$ has $S_1$-type. 
By complex interpolation with respect to $L_2(\mu, S^n_2)$ we obtain that $L_p(\mu, S^n_p)$($1 \leq p \leq 2$) has $S_p$-type.

Now we consider the case $2\leq p < \infty$. We can show that $L_p(\mu, S^n_p)$ has $S_2$-type by the direct calculation as
above. Since $\G_2 \stackrel{\F_\R}{\longrightarrow} S_2 \stackrel{I_{2,\infty}}{\longrightarrow} S_{\infty}$ is a
contraction, where $I_{2,\infty}$ is the corresponding formal identity, so is
$$\F_\R \otimes I_{L_{\infty}(\mu, S^n_{\infty})} : \G_2\otimes_{\lambda} L_{\infty}(\mu, S^n_{\infty})
\rightarrow S_{\infty}\otimes_{\lambda}L_{\infty}(\mu, S^n_{\infty}),$$ where $\otimes_{\lambda}$ is the injective tensor
product in the category of Banach spaces. Note that we have the following contraction
$$\G_2(E) \subseteq L_2(\Om,E) \hookrightarrow L_2\Om\otimes_{\lambda} E \stackrel{P_2\otimes I_E
}{\longrightarrow} \G_2\otimes_{\lambda} E,$$ for $E = L_{\infty}(\mu, S^n_{\infty})$ and $P_2$ is the gaussian
projection from $L_2(\Om)$ onto $\G_2$. Since $$B(S_1(L_1(\mu)), S^n_{\infty}) \supseteq S_{\infty}\otimes_{\lambda}L_{\infty}(\mu,
S^n_{\infty}) \cong S_{\infty}(L_{\infty}(\mu, S^n_{\infty})) \subseteq CB(S_1(L_1(\mu)), S^n_{\infty})$$ isomorphically, 
we have a bounded map $$\F_\R \otimes I_{L_{\infty}(\mu, S^n_{\infty})} : \G_2(L_{\infty}(\mu, S^n_{\infty}))
\rightarrow S_{\infty}(L_{\infty}(\mu, S^n_{\infty})),$$ which implies $L_{\infty}(\mu, S^n_{\infty})$ has $S_{\infty}$-cotype. 
By complex interpolation with respect to $L_2(\mu, S^n_2)$ we obtain that $L_p(\mu, S^n_p)$($2 \leq p \leq \infty$) has $S_p$-cotype.

The other statements concerning best $S_p$-type and $S_q$-cotype follows by Remark \ref{rem-def} and the Banach space case, 
except the fact that infinite dimensional $L_{\infty}(\mu)$ does not have $S_1$-type. For simplicity we just consider the case $c_0$, 
the space of null sequences. Note that by the dominance of the gaussian average over the Rademacher average (for example, (4.2) of \cite{TJ}) we have
\begin{align*}
\begin{split}
T^{o,n}_1(c_0) & = \norm{S^n_1(c_0) \rightarrow L_2(\Om, c_0)\,\, , \,\, (x_{ij}) \mapsto \sum^n_{i,j=1}g_{ij}\otimes x_{ij}}\\
& \gtrsim \norm{S^n_1(c_0) \rightarrow L_2([0,1], c_0)\,\, , \,\, (x_{ij}) \mapsto \sum^n_{i,j=1}r_{ij}\otimes x_{ij}}\\
& = \norm{L_2([0,1], \ell_1) \rightarrow S^n_{\infty}(\ell_1)\,\,
, \,\, f \mapsto (\left\langle r_{ij}, f \right\rangle)_{i,j}},
\end{split}
\end{align*}
where $\{r_{ij}\}$ is an re-indexing of the classical Rademacher system $\{ r_i \}$.

Set $f(t) = \prod^n_{i,j=1}(1 + r_{ij}(t)r_{ij}) \in L_2([0,1],L_1[0,1])$.
Then we have $$\norm{f(t)}_{L_1[0,1]} = \int^1_0 \abs{\prod^n_{i,j=1}(1 + r_{ij}(t)r_{ij}(s))}ds
= \int^1_0\prod^n_{i,j=1}(1 + r_{ij}(t)r_{ij}(s))ds = 1$$ for all $t \in [0,1]$, and consequently $$\norm{f}_{L_2([0,1],L_1[0,1])} = 1.$$
On the other hand, we have $\left\langle r_{ij}, f \right\rangle = r_{ij}$, so that
\begin{align*}
\begin{split}
\norm{(\left\langle r_{ij}, f
\right\rangle)_{i,j}}_{S^n_{\infty}(L_1[0,1])} & =
\norm{\sum^n_{i,j=1}e_{ij}\otimes r_{ij}}_{S^n_{\infty}(L_1[0,1])}
 \sim \norm{\sum^n_{i,j=1}e_{ij}\otimes e_{ij}}_{S^n_{\infty}(R_{n^2}+C_{n^2})}\\
& = \norm{R_{n^2}\cap C_{n^2} \rightarrow S^n_{\infty} \,\, , \,\,
\delta_{ij} \mapsto e_{ij}}_{cb}\\
& \geq \frac{\norm{\sum^n_{i,j=1}e_{ij}\otimes
e_{ij}}_{S^n_{\infty}(S^n_{\infty})}}{\max\Big\{
\norm{\sum^n_{i,j=1}e_{ij}e^*_{ij}}^{\frac{1}{2}}_{\min},
\norm{\sum^n_{i,j=1}e^*_{ij}e_{ij}}^{\frac{1}{2}}_{\min}\Big\}}\\
& = n^{\frac{1}{2}}.
\end{split}
\end{align*}
Since $\text{span}\{r_A = \prod_{(i,j)\in A}r_{ij} : A \subseteq \{(i,j) | 1\leq i,j \leq n\} \}$ in $L_1[0,1]$ is completely
isometric to $\ell^{2^{n^2}}_1$ we get $$T^{o,n}_1(c_0) \gtrsim n^{\frac{1}{2}},$$ the desired result.

\end{proof}

\begin{rem}\label{rem-trivial}{ \rm
We do not need $\sigma$-finiteness of $\mu$ to prove $L_1(\mu, S^n_1)$ has $S_1$-type and $L_{\infty}(\mu, S^n_{\infty})$ has
$S_{\infty}$-cotype in the above theorem, and we can similarly show that every maximal operator space has $S_1$-type and every minimal
operator space has $S_{\infty}$-cotype.}
\end{rem}

Now, we consider $S_p$-type and $S_q$-cotype of infinite dimensional Schatten classes. 
Unfortunately we could not determine the best $S_p$-type and $S_q$-cotype of those spaces at the time of this writing,
but the following estimate shows that they are quite different from type and cotype as Banach spaces.

\begin{lem} \label{lem-Sptype-estimate-Sp}
Let $1\leq p \leq 2$. Then we have $$T_{S_p,n}(S_p) \sim C_{S_{p'},n}(S_{p'}) \sim n^{\frac{1}{p} - \frac{1}{2}}.$$
\end{lem}
\begin{proof}
First, we consider type constants. By Theorem 9.8.5 of \cite{P5} we have
\begin{align*}
\begin{split}
T_{S_p,n}(S_p) & \sim \norm{\F^{-1}_{\R} \otimes I_{S_p} : S^n_p(S_p) \rightarrow \G^n_p(S_p)}
= \norm{S^n_p \rightarrow \G^n_p , \,\, e_{ij} \mapsto r_{ij}}_{cb}\\
& \sim \norm{R_{n^2}[p'] \cap C_{n^2}[p'] \rightarrow S^n_{p'} ,\,\, e_{ij} \mapsto e_{ij}}_{cb}
\leq \norm{R_{n^2}[p'] \rightarrow S^n_{p'} , \,\, e_{ij} \mapsto e_{ij}}_{cb}\\
& = \norm{R_n[p'] \otimes_h R_n[p']\rightarrow C_n[p'] \otimes_h R_n[p'] , \,\, e_{1i}\otimes e_{1j} \mapsto e_{i1}\otimes e_{1j}}_{cb}\\
& \leq \norm{R_n[p']\rightarrow C_n[p'] , \,\, e_{1i} \mapsto e_{i1}}_{cb}.
\end{split}
\end{align*}
By complex interpolation we get for $\frac{\theta}{2} + \frac{1-\theta}{\infty} = \frac{1}{p'}$
\begin{align*}
\begin{split}
\lefteqn{\norm{R_n[p']\rightarrow C_n[p'] , \,\, e_{1i} \mapsto
e_{i1}}_{cb}}\\ & \leq \norm{R^n_{2}\rightarrow C^n_{2} , \,\, e_{1i}
\mapsto e_{i1}}^{\theta}_{cb}
\norm{R^n_{\infty}\rightarrow C^n_{\infty} , \,\, e_{1i} \mapsto e_{i1}}^{1-\theta}_{cb}
= n^{\frac{1}{p} - \frac{1}{2}}.
\end{split}
\end{align*}
For the upper bound we consider
\begin{align*}
\begin{split}
\lefteqn{\norm{R_{n^2}[p'] \cap C_{n^2}[p'] \rightarrow S^n_{p'} , \,\, e_{ij} \mapsto e_{ij}}_{cb}}\\
& \geq \frac{\norm{\sum^n_{i,j=1}e_{ij}\otimes e_{ij}}_{S^n_{p'}(S^n_{p'})}}{\max \Big\{ \norm{(\sum^n_{i,j=1}e^*_{ij} e_{ij})^{\frac{1}{2}}}_{S^n_{p'}},\norm{(\sum^n_{i,j=1}e_{ij} e^*_{ij})^{\frac{1}{2}}}_{S^n_{p'}} \Big\}}
= n^{\frac{1}{p} - \frac{1}{2}}
\end{split}
\end{align*}

Now we consider cotype constants. For $1< p \leq 2$ we have
\begin{align*}
\begin{split}
C_{S_{p'},n}(S_{p'}) & \sim \norm{\F_{\R} \otimes I_{S_p} : \G^n_{p'}(S_{p'}) \rightarrow S^n_{p'}(S_{p'})}\\
& \sim \norm{R_{n^2}[p'] \cap C_{n^2}[p'] \rightarrow S^n_{p'} ,
\,\, e_{ij} \mapsto e_{ij}}_{cb} = n^{\frac{1}{p} - \frac{1}{2}}.
\end{split}
\end{align*}
For $p = 1$ we have by Proposition 45.1 in \cite{TJ}
\begin{align*}
\begin{split}
C^{o,n}_{\infty}(S_{\infty}) & \sim \norm{\F_{\R} \otimes I_{S_{\infty}} : \G^n_1(S_{\infty}) \rightarrow S^n_{\infty}(S_{\infty})}\\
& \geq \frac{\norm{\sum^n_{i,j=1} e_{ij}\otimes e_{ij}}_{S^n_{\infty}(S^n_{\infty})}}{\int_{\Om} \norm{(g_{ij}(\om))}_{S^n_{\infty}}dP(\om)}
\gtrsim n^{\frac{1}{2}}.
\end{split}
\end{align*}
We get the upper bound by (3) of Remark \ref{rem-def}.

\end{proof}

\begin{thm}

Let $\frac{1}{p} + \frac{1}{p'} = 1$.

\begin{itemize}
\item[(1)]
When $1 \leq p < \frac{4}{3}$, $S_p$ does not have $S_1$-type nor $S_s$-cotype for $2\leq s < p'$.

\item[(2)]
When $\frac{4}{3} \leq p < 2$, $S_p$ does not have $S_r$-type nor $S_s$-cotype for $$\frac{2p}{4-p} < r \leq 2\leq s < p'.$$

\item[(3)]
When $2 < p \leq 4$, $S_p$ does not have $S_r$-type for  nor $S_s$-cotype for $$p' < r \leq 2 \leq s < \frac{2p}{4-p}.$$

\item[(4)]
When $4 < p < \infty$, $S_p$ does not have $S_r$-type for $p' < r \leq 2$ nor $S_{\infty}$-cotype.

\item[(5)] $S_{\infty}$ does not have $S_1$-type nor $S_{\infty}$-cotype.

\end{itemize}
\end{thm}
\begin{proof}
Since the formal identity $S^n_p(E) \rightarrow S^n_r(E)$ has norm $\leq n^{\frac{1}{r} - \frac{1}{p}}$ for $1\leq r < p$, 
we get $$T^{o,n}_r(S_p) \geq T_{S_p,n}(S_p)n^{\frac{1}{p} - \frac{1}{r}} \gtrsim n^{\frac{2}{p} - \frac{1}{2} - \frac{1}{r}},$$
which means $S_p$ does not have $S_r$-type for $\frac{2}{p} - \frac{1}{2} - \frac{1}{r} > 0 \Leftrightarrow \frac{2p}{4-p}<r$.
The other statements are obtained by duality (Proposition \ref{duality}), Theorem \ref{Thm-SemiNoncomLp},
Theorem \ref{compute1} and the fact that $R[p]$ is a closed subspace of $S_p$.
\end{proof}

We close this section with the case of $C^*$-algebras and their duals. 
$S_1$-type and $S_{\infty}$-cotype are related to subhomogeneity of a $C^*$-algebra.
\begin{thm}\label{thm-subhom}
Let $A$ be a $C^*$-algebra. Then $A$ is subhomogeneous if and only if $A$ has $S_{\infty}$-cotype if and only if $A^*$ has $S_1$-type.
Moreover, if $A$ is not subhomogeneous, then we have $$T^{o,n}_{1}(A^*) \sim C^{o,n}_{\infty}(A) \sim n^{\frac{1}{2}},$$ which is worst possible.
\end{thm}
\begin{proof}
Suppose that $A$ is $m$-subhomogeneous for some $m \in \mathbb{N}$. Then we can assume that $A^{**}$ is a subalgebra of
$L_{\infty}(\mu, S^m_{\infty})$ for some measure space $(\mathcal{M}, \mu)$. Since $L_{\infty}(\mu, S^m_{\infty})$ has
$S_{\infty}$-cotype and $L_1(\mu, S^m_1)$ has $S_1$-type (Remark \ref{rem-trivial}), so does $A$ and $A^*$, respectively.

Now we assume that $A$ is not subhomogeneous, then for any $\epsilon>0$ and $n \geq 1$ there are completely positive maps
$$\rho : S^n_{\infty} \rightarrow A\,\, \text{and}\,\, \sigma : A \rightarrow S^n_{\infty}\,\, \text{such that}\,\, 
\norm{\sigma \rho - I_{S^n_{\infty}}}_{cb} \leq n\cdot\norm{\sigma \rho - I_{S^n_{\infty}}}< \frac{\epsilon}{n}$$
by Lemma 2.7 of \cite{S} and \cite{HT}. Then we have
$$\int_{\Om} \norm{\sum^n_{i,j=1}g_{ij}(\om)\rho e_{ij}}_{A}dP(\om) \leq \norm{\rho}_{cb} \int_{\Om}
\norm{\sum^n_{i,j=1}g_{ij}(\om)e_{ij}}_{S^n_{\infty}}dP(\om) \lesssim n^{\frac{1}{2}}$$ and
\begin{align*}
\begin{split}
\norm{(\rho e_{ij})}_{S^n_{\infty}(A)} & \geq \norm{\sigma}^{-1}_{cb}\norm{(\sigma\rho e_{ij})}_{S^n_{\infty}(S^n_{\infty})}\\
& = \norm{((\sigma \rho - I_{S^n_{\infty}})e_{ij} + e_{ij})}_{S^n_{\infty}(S^n_{\infty})}\\
& \geq \norm{(e_{ij})}_{S^n_{\infty}(S^n_{\infty})} - \norm{((\sigma \rho - I_{S^n_{\infty}})e_{ij})}_{S^n_{\infty}(S^n_{\infty})}\\
& \geq n - \norm{\sigma \rho - I_{S^n_{\infty}}}_{cb}\cdot n \geq n - \epsilon,
\end{split}
\end{align*}
which implies $C^{o,n}_{\infty}(A) \sim n^{\frac{1}{2}}$. $T^{o,n}_{1}(A^*) \sim n^{\frac{1}{2}}$ is obtained by duality.

\end{proof}

\section{Type $(p,H)$ and cotype $(q,H)$ of operator spaces}\label{sec-type(p,H)-cotype(q,H)}

\subsection{Definitions and basic properties}

We fix a subquadratic homogeneous Hilbertian operator space $H$ from now on.
Now for an operator space $E$ we define $(2,H)$-summing norm of a map $v : E\rightarrow \ell_2$ by
$$\pi_{2,H}(v) = \sup\Big\{ \frac{(\sum_k \norm{vx_k}^2 )^{\frac{1}{2}}}{\norm{\sum_k x_k\otimes e_k}_{E\otimes_{\text{min}}H}} \Big\}.$$
Note that the subquadraticity of $H$ ensure that $\pi_{2,H}(\cdot)$ is actually a norm. (p.82 of \cite{P2})
Also note that all results remain true for all $H$ which is completely isomorphic to a subquadratic homogeneous Hilbertian operator space
allowing suitable constants.

\begin{defn}
An operator space $E$ is called type $(2,H)$ if there is a constant $C>0$ such that 
$$\ell(u) \leq C\pi_{2,H}(u^*)$$ for all $n\in \N$ and $u : \ell^n_2 \rightarrow E$.

$E$ is called cotype $(2,H)$ is there is a constant $C'>0$ such that $$\ell^*(v) \leq C'\pi_{2,H}(v)$$
for all $n\in \N$ and $v : E\rightarrow \ell^n_2$. We denote the infimums of such $C$ and $C'$ by $T_{2,H}(E)$ and $C_{2,H}(E)$, respectively.
\end{defn}

We give a description of the trace dual of $\pi_{2,H}$.
\begin{prop}
For $n \in \mathbb{N}$ and $u : \ell^n_2 \rightarrow E$ we have $$\pi^*_{2,H}(u) = \inf \norm{A}_{HS}\norm{B}_{cb},$$ where the infimum runs over all
possible factorization $u : \ell^n_2 \stackrel{A}{\longrightarrow} H^* \stackrel{B}{\longrightarrow} E$.

\end{prop}\label{prop-tracedual-pi2H}
\begin{proof}
Let $\alpha(u)$ be the infimum on the right hand side. For any $v : E \rightarrow \ell^n_2$ and factorization 
$u : \ell^n_2 \stackrel{A}{\longrightarrow} H^* \stackrel{B}{\longrightarrow} E$ we have
\begin{align*}
\begin{split}
\abs{\text{tr}(vu)} & \leq \abs{\text{tr}(vBA)} = \abs{\text{tr}(AvB)} \leq \norm{A}_{HS}\norm{vB}_{HS}\\
& = \norm{A}_{HS}\pi_{2,H}(vB) \leq \norm{A}_{HS}\norm{B}_{cb}\pi_{2,H}(v),
\end{split}
\end{align*}
which implies $\pi^*_{2,H}(u) \leq \alpha(u)$.

For the converse inequality we will show that $\alpha^*(v) \geq \pi_{2,H}(v)$ for any $v : E \rightarrow \ell^n_2$.
For any given $\epsilon > 0$ we choose $(x_k) \subseteq E$ such that 
$$\Big(\sum_k \norm{vx_k}^2\Big)^{\frac{1}{2}} \geq (1-\epsilon) \norm{\sum_k x_k\otimes e_k}_{E\otimes_{\text{min}}H}.$$
Then there is $(y^*_k) \in \ell^n_2(\ell^n_2)$ with norm 1 such that 
$$\abs{\left\langle (y^*_k), (vx_k)\right\rangle} \geq (1-\epsilon)^2 \norm{(vx_k)}_{\ell^n_2(\ell^n_2)} 
= (1-\epsilon)^2\Big(\sum_k \norm{vx_k}^2\Big)^{\frac{1}{2}}.$$
We set $A : \ell^n_2 \rightarrow H^*_n,\; z\mapsto (\left\langle y^*_k, z\right\rangle)_k$ and 
$B: H^*_n \rightarrow E,\; e_k \mapsto x_k$. Then we have
\begin{align*}
\begin{split}
\alpha^*(v) & \geq \frac{\abs{\text{tr}(vBA)}}{\norm{A}_{HS}\norm{B}_{cb}} = 
\frac{\abs{\text{tr}(AvB)}}{\norm{\sum_k x_k\otimes e_k}_{E\otimes_{\text{min}}H}}\\
& = \frac{\abs{\sum_k \left\langle A^* e_k, vBe_k \right\rangle}}{\norm{\sum_k x_k\otimes e_k}_{E\otimes_{\text{min}}H}}
= \frac{\abs{\left\langle (y^*_k), (vx_k)\right\rangle}}{\norm{\sum_k x_k\otimes e_k}_{E\otimes_{\text{min}}H}}\\
& \geq (1-\epsilon)^2\pi_{2,H}(v).
\end{split}
\end{align*}

\end{proof}

Now we consider the k-th c.b. approximation number of $T : E \rightarrow F$ by
$$a^o_k(T) := \inf \{ \norm{T-S}_{cb} : S \in CB(E,F), \text{rk}(S) < k\}.$$
(See \cite{O} for operator space versions of Gelfand and Kolmogorov
numbers.)

\begin{prop}
For $u : H^*_n \rightarrow E$ we have $$\Big(\sum_k a^o_k(u)^2 \Big)^{\frac{1}{2}} \leq \pi^*_{2,H}(u).$$
\end{prop}
\begin{proof}
By Proposition \ref{prop-tracedual-pi2H} for any given $\epsilon > 0$ we have a factorization
$$u : H^*_n \stackrel{A}{\longrightarrow} H^* \stackrel{B}{\longrightarrow} E$$ with 
$\norm{A}_{HS}\norm{B}_{cb} \leq (1+\epsilon)\pi^*_{2,H}(u).$ Thus, we have
\begin{align*}
\begin{split}
\Big(\sum_k a^o_k(u)^2 \Big)^{\frac{1}{2}} & = \Big(\sum_k a^o_k(BA)^2 \Big)^{\frac{1}{2}} 
\leq \norm{B}_{cb}\Big(\sum_k a^o_k(A)^2 \Big)^{\frac{1}{2}}\\ & = \norm{B}_{cb}\norm{A}_{HS} \leq (1+\epsilon)\pi^*_{2,H}(u).
\end{split}
\end{align*}

\end{proof}

Recall that there is a constant $C>0$ such that 
$$\Big(\sum_k a_k(u)^2 \Big)^{\frac{1}{2}} \leq \pi_{q,2}(u) \leq \frac{Cq}{q-2}\Big(\sum_k a_k(u)^2 \Big)^{\frac{1}{2}}$$
for any $u : \ell^n_2 \rightarrow X$ $n \in \mathbb{N}$ and $2<q \leq \infty$. This equivalence and (\ref{def-type-p}) and (\ref{def-cotype-q}) 
lead us to the following definition.
\begin{defn} Let $1\leq p < 2 < q \leq \infty$.
An operator space $E$ is called type $(p,H)$ if there is a constant $C>0$ such that 
$$\Big(\sum_k a^o_k(v)^{p'} \Big)^{\frac{1}{p'}} \leq C \cdot \ell^*(v)$$ for all $n\in \N$ and $v : E\rightarrow H_n$.

$E$ is called cotype $(q,H)$ is there is a constant $C'>0$ such that
$$\Big(\sum_k a^o_k(u)^q \Big)^{\frac{1}{q}} \leq C' \cdot \ell(u)$$ for all $n\in \N$ and $u : H^*_n \rightarrow E$.
We denote the infimums of such $C$ and $C'$ by $T_{p,H}(E)$ and $C_{q,H}(E)$, respectively.
\end{defn}

\begin{rem}
{\rm 
\begin{itemize}
\item[(1)]
It is clear from the definition that type $(p,H)$ and cotype $(q,H)$ imply type $p$ and cotype $q$ as Banach spaces, respectively.
\item[(2)]
Let $$S^o_r(E,F) = \{ u \in CB(E,F) : \norm{(a^o_k(u))_{k\geq 1}}_{\ell_r} < \infty \}$$ for $1\leq r \leq \infty$.
Then by the same argument in the proof of Proposition 1 in \cite{KTz} we have 
$$K(t, u ; S^o_1(E,F), S^o_{\infty}(E,F)) \sim K(t, (a^o_k(u))_{k\geq 1} ; \ell_1, \ell_{\infty})$$
for $t > 0$, where $K(t, \cdot ;E_0, E_1)$ implies the $K$-functional with respect to a compatible pair of Banach spaces $(E_0,E_1)$. 
Thus we have  $$[S^o_2(E,F), S^o_{\infty}(E,F)]_{\frac{2}{q}} = S^o_q(E,F)$$ for $2 < q < \infty$.
When $E$ is cotype $(2,H)$ we have $$\Big(\sum_k a^o_k(u)^2 \Big)^{\frac{1}{2}} \leq C_{2,H}(E) \ell(u)$$ 
for all $n\in \N$ and $u : H^*_n \rightarrow E$, thus cotype $(q,H)$ (resp. type $(p,H)$) behaves well via interpolation 
as in the Banach space case.
\end{itemize}
}
\end{rem}

As in $S_p$-type and $S_q$-cotype case we have the following duality results.

\begin{prop} \label{prop-duality-(p,H)} Let $E$ be an operator space and $1\leq p \leq 2$.
\begin{itemize}
\item[(1)] If $E$ has type $(p,H)$, then $E^*$ has cotype $(p', H)$ with $$C_{p', H}(E^*) \leq T_{p,H}(E).$$ 

\item[(2)] If $E$ has cotype $(p',H)$ and is $K$-convex as a Banach space, then $E^*$ has type $(p, H)$ with 
$$T_{p,H}(E^*)\leq K(E)C_{p',H}(E).$$
\end{itemize}
\end{prop}
\begin{proof}
Note that we have $\ell^*(v) \leq \ell(v^*)$ and $\ell(v^*) \leq K(X)\ell^*(v)$ for any Banach space $X$ and $v : \ell_2 \rightarrow X$.
\end{proof}

\subsection{The case of homogeneous Hilbertian operator spaces}

If we consider type $(p,H)$ and cotype $(q,H)$ of homogeneous Hilbertian operator spaces, then the calculation becomes simple,
so that we can completely determine type and cotype in some cases. We only consider cotype cases, since type cases can be directly obtained by duality.

Let's start with the following lemma about the approximation number of formal identities between homogeneous Hilbertian operator spaces.
Recall that the $k$-th c.b. Gelfand number of $u :E \rightarrow F$ between operator spaces is defined by
$$c^o_k(u) := \inf \{ \norm{u|_S}_{cb} : S \subseteq E,\; \text{codim}S < k\}$$ for $k\in \mathbb{N}$, and clearly we have 
\begin{equation}\label{GelfandLessApprox}
c^o_k(u) \leq a^o_k(u).
\end{equation}

\begin{lem} \label{lem-diagonal}
Let $\Hi$ and $\Hi'$ be homogeneous Hilbertian operator spaces. Then for the $n$-dimensional formal identity $id_n : \Hi_n \rightarrow \Hi'_n$
we have $$c^o_k(id_n) = a^o_k(id_n) = \norm{id_{n-k+1} : \Hi_{n-k+1} \rightarrow \Hi'_{n-k+1}}_{cb}$$ for $1\leq k \leq n$.
\end{lem}
\begin{proof}
Fix $1\leq k \leq n$. By (\ref{GelfandLessApprox}) it is enough to show that 
$$c^o_k(id_n)\geq \norm{id_{n-k+1} : \Hi_{ n-k+1} \rightarrow \Hi'_{n-k+1}}_{cb}.$$
Now we consider any subspace $S \subseteq \Hi_{n}$ with $m := \text{dim}S = n - \text{codim}S \geq n-k+1$. 
Then their is a partial isometry $U_S : \ell^m_2 \rightarrow \ell^n_2$ whose image is $S$ and $$U^*_S U_S = I_{\ell^m_2}.$$
Thus, we have 
\begin{align*}
\begin{split}
\norm{id_{n-k+1}}_{cb} & \leq \norm{id_m : \Hi_{ m} \rightarrow \Hi'_{ m}}_{cb}\\
& \leq \norm{U_S : \Hi_{m} \rightarrow \Hi'_{n}}_{cb} \norm{U^*_S : \Hi'_{n} \rightarrow \Hi'_{m}}_{cb}\\
& = \norm{id_n|_S :\Hi_{n} \rightarrow \Hi'_{n}}_{cb}.
\end{split}
\end{align*}

\end{proof}

\begin{prop}\label{prop-cotype(2,H)-Hilbertian}
Let $\Hi$ be a homogeneous Hilbertian operator spaces. Then, $\Hi$ is cotype $(2,H)$ if and only if the formal identity
$$id : H^* \rightarrow \Hi$$ is completely bounded.
\end{prop}
\begin{proof}
The necessity part is clear from the definition. Now we suppose that $\Hi$ is cotype $(2,H)$.
Then for $n \in \mathbb{N}$ and $u : H^*_n\rightarrow \Hi_n$ we have 
$$\sum_k a^o_k(u)^2 \leq C^2 \ell^2(u),$$ where $C = C_{2,H}(\Hi)$.
If we set $u = id_n : H^*_n\rightarrow \Hi_n$, then by Lemma \ref{lem-diagonal} we have 
$$ C^2 n \geq \sum^n_{k=1} \norm{id_k : H^*_k \rightarrow \Hi_k}_{cb}^2 \geq \sum^n_{k=[\frac{n}{2}]} \norm{id_k}_{cb}^2
\geq \frac{n}{2}\norm{id_{[\frac{n}{2}]}}_{cb}^2,$$ which means $\norm{id_{[\frac{n}{2}]}}_{cb} \leq \sqrt{2}C$
and consequently $\norm{id : H^* \rightarrow \Hi}_{cb}$ is bounded.
\end{proof}

Now we focus on $R[p]$ and $C[p]$ case.
\begin{thm}
Let $1\leq p, q \leq \infty$. Then $R[q]$ (resp. $C[q]$) is cotype $(s,C[p])$ (resp. $(s,R[p])$) 
if and only if $\abs{\frac{1}{p} -\frac{1}{q}} + \frac{1}{s} \leq \frac{1}{2}.$ 
\end{thm}
\begin{proof}
Consider $u : R[p] \rightarrow R[q]$. Let $\abs{\frac{1}{p} -\frac{1}{q}} = \frac{1}{r}$, then by Lemma 5.9 of \cite{Xu} 
we have $$CB(R[p], R[q]) \cong S_r$$ isometrically. Since c.b. approximation numbers of $u$ and $\ell(u)$ are both unitarily invariant 
we can assume that $$u = \text{diag}(u_1, u_2, \cdots, u_n), \; n\in \mathbb{N}$$ with
$\abs{u_1} \geq \abs{u_2} \geq \cdots \geq \abs{u_n}$ by a usual density argument. 

Now we suppose $R[q]$ is cotype $(2,C[p])$ and set $u_1 = \cdots = u_n = 1$. Then we have 
$$\ell(u) = \Big(\sum^n_{k=1}\abs{u_k}^2\Big)^{\frac{1}{2}} = n^{\frac{1}{2}}$$ and by Lemma \ref{lem-diagonal}
$$\Big(\sum_k a^o_k(u)^s \Big)^{\frac{1}{s}} = \Big(\sum_k (n-k+1)^{\frac{s}{r}} \Big)^{\frac{1}{s}} \sim n^{\frac{1}{r} + \frac{1}{s}}.$$
Consequently, we have $$\abs{\frac{1}{p} -\frac{1}{q}} + \frac{1}{s} = \frac{1}{r} + \frac{1}{s} \leq \frac{1}{2}.$$

For the converse we observe the following.
$$\Big(\sum_k a^o_k(u)^s \Big)^{\frac{1}{s}} \leq \Big[\sum_k \Big(\sum_{i\geq k}\abs{u_i}^r \Big)^{\frac{s}{r}}\Big]^{\frac{1}{s}}
= \norm{U}_{\ell^n_s(\ell^n_r)},$$ where $U=(u_{ij})^n_{i,j=1}$ with $u_{ij} = u_j$ for $j\leq i$ and $u_{ij} = 0$ elsewhere.
Thus, it is enough to show that $$\norm{U}_{\ell^n_s(\ell^n_r)} \leq \Big(\sum^n_{k=1}\abs{u_k}^2\Big)^{\frac{1}{2}}$$
for $\frac{1}{r} + \frac{1}{s} = \frac{1}{2}$. Since we have 
$$\ell^n_s(\ell^n_r) = [\ell^n_4(\ell^n_4), \ell^n_2(\ell^n_{\infty})]_{\theta}$$ for $r \geq s$ and $\theta = 1-\frac{4}{r}$ 
and $$\ell^n_s(\ell^n_r) = [\ell^n_4(\ell^n_4), \ell^n_{\infty}(\ell^n_2)]_{\psi}$$ for $r < s$ and $\psi = 1-\frac{4}{s}$ it suffices to consider
the following three extremal cases : $(r,s) = (2, \infty), (\infty, 2)$ and $(4,4)$.

When $(r,s) = (2, \infty)$ or $(\infty, 2)$ it is trivial from the definition. The case $r = s = 4$ is obtained from the following.
$$\Big(\sum^n_{k=1} k \abs{u_k}^4 \Big)^{\frac{1}{4}} \leq \Big(\sum^n_{k=1} \abs{u_k}^2 \Big)^{\frac{1}{2}}.$$
Indeed, we can show the above inequality by induction on $n$. 
When $n=1$ it is trivial. Suppose that it is true for $n$, then we have
$$\sum^n_{k=1} k \abs{u_k}^4 \leq \Big(\sum^n_{k=1} \abs{u_k}^2 \Big)^2,$$ and consequently
\begin{align*}
\begin{split}
\sum^{n+1}_{k=1} k \abs{u_k}^4 & = \sum^n_{k=1} k \abs{u_k}^4 + (n+1)\abs{u_{k+1}}^4 
\leq \sum^n_{k, l=1} \abs{u_k}^2\abs{u_l}^2 +(n+1)\abs{u_{k+1}}^4 \\ & \leq \sum^{n+1}_{k, l=1} \abs{u_k}^2\abs{u_l}^2
= \Big(\sum^{n+1}_{k=1} \abs{u_k}^2 \Big)^2,
\end{split}
\end{align*}
since $\abs{u_k}$'s are non-increasing.

The proof for the $C[q]$ is the same.  
\end{proof}

\begin{rem}
{\rm 
Since cotype $(2,H)$ is a local property if $F$ is $\lambda$-c.b. representable in $E$ for some $\lambda > 0$
(i.e. every finite dimensional subspace of $F$ can be $(1+\epsilon)\lambda$-c.b. embedded in $E$ for any $\epsilon > 0$) then
cotype $(2,H)$ property of $E$ can be transferred to $F$. However, sometimes cotype $(2,H)$ property can be transferred to an operator spaces
related in a weaker sense. More precisely, let's say that ``$F$ is $\lambda$-representable in $E$ at every matrix level" i.e. 
for any $m\in \N$, $\epsilon > 0$ and finite dimensional subspace $F'$ of $F$
there is a subspace $E' \subseteq E$ and an isomorphism $T : F' \rightarrow E'$ such that
$$\norm{I_{M_m}\otimes T^{-1} : M_m(E') \rightarrow M_m(F')} = 1$$ and $$\norm{I_{M_m}\otimes T : M_m(F') \rightarrow M_m(E')}\leq (1+ \epsilon)\lambda.$$
Then cotype $(2,H)$ property of $E$ can be transferred to $F$ if $F$ is another homogeneous Hilbertian operator space.

Indeed, by Proposition \ref{prop-cotype(2,H)-Hilbertian} we need to check that $\norm{id_n : H^*_n \rightarrow F_n}_{cb}$
is uniformly bounded with respect to $n \in \N$. Now we fix $n \in \N$. Then for any $\epsilon > 0$ there is $m \in \N$ such that
$$\norm{id_n : H^*_n \rightarrow F_n}_{cb} \leq (1+ \epsilon)\norm{I_{M_n}\otimes id_n : M_m(H^*_n) \rightarrow M_m(F_n)}.$$
Now we set $F' = F_{2n}$ and choose $E' \subseteq E$ and $T$ as above.
Since $E'$ is a subspace of $E$ we have $$\sum^{2n}_{k=1}a^o_k(u) \leq C_{2,H}(E)\ell(u)$$ for any $u : H^*_{2n} \rightarrow E'$.
For $u = T\circ id_{2n}$ we have $$\ell(u) \leq (1+\epsilon)\lambda\ell(id_{2n}) = (1+\epsilon)\lambda \sqrt{2n}$$ and
$$\sum^{2n}_{k=1}a^o_k(u) \geq \sum^{2n}_{k=1}a^{M_m}_k(u) \geq \sum^{2n}_{k=1}a^{M_m}_k(id_{2n}),$$ where
\begin{align*}
\begin{split}
\lefteqn{a^{M_m}_k(v : E_1 \rightarrow E_2)}\\ 
& := \inf \{ \norm{I_{M_m}\otimes(v-w) : M_m(E_1) \rightarrow M_m(E_2)} : w \in B(E_1,E_2), \text{rk}(w) < k\}.
\end{split}
\end{align*}
By a similar argument as in the Lemma \ref{lem-diagonal} we get 
$$a^{M_m}_k(id_{2n}) = \norm{I_{M_m}\otimes id_k : M_m(H^*_{2n})\rightarrow M_m(F_{2n})}$$ and consequently
$$\norm{I_{M_m}\otimes id_n : M_m(H^*_n) \rightarrow M_m(F_n)}
\leq \sqrt{2}(1+\epsilon)\lambda.$$

The situation as above happens. By \cite{P7} we know that for every infinite dimensional operator space $E$ there is a homogeneous Hilbertian operator space contained in $E^\U$, an ultrapower of $E$.
It is well known that the local structure of $E^\U$ as an operator space is not the same as $E$ unlike in Banach space case.
However, by a similar argument as in the Banach space case we can show that $E^\U$ is``$\lambda$-representable in $E$ at every matrix level" for $\lambda = 1$.
}
\end{rem}

\subsection{The case of $L_p$ spaces}

As in the Banach space case type $(1,H)$ and cotype $(\infty, H)$ are trivial for certain $H$.
\begin{prop}\label{prop-type(1,H)}
Every operator space has type $(1,H)$ and cotype $(\infty, H)$ for $H = R[p], C[p]$ and $RC[p]$.
\end{prop}
\begin{proof}
We only prove the type case since the cotype case is obtained by duality.
Note that an operator space $E$ is type $(1,H)$ if and only if there is a constant $C>0$ such that $$\norm{v}_{cb} \leq C\cdot\ell^*(v)$$
for any $n \in \mathbb{N}$ and $v : E \rightarrow H_n$.

First, we consider the case $H = R$. Since we have $$\norm{v : E \rightarrow R_n}_{cb} \leq \norm{v : \min(E) \rightarrow R_n}_{cb} 
= \pi_2(v : E\rightarrow \ell^n_2)$$ by (1.45) of \cite{BLe} it is enough to show that $$\pi_2(v) \leq \ell^*(v)$$ for any $n \in \mathbb{N}$ and 
$v : E \rightarrow \ell^n_2$. By trace duality it is equivalent to the following well-known result (for example (3.14) of \cite{P1})
$$\ell(u) \leq \pi_2(u)$$ for any $n \in \mathbb{N}$ and $u : \ell^n_2 \rightarrow E$.

We can prove the case $H = C$ with the same proof and by combining these two result we get the case $H = R\cap C$ and $R+C$. 
Finally, we are done by interpolation.
\end{proof}

Using Proposition \ref{prop-cotype(2,H)-Hilbertian} we can determine the condition for $L_p$ ($1\leq p \leq 2$) spaces to be cotype 2.

\begin{thm}\label{thm-cotype(2,H)-Lp}
Let $1\leq p \leq 2$ and $\mu$ be a $\sigma$-finite measure.
\begin{itemize}
\item[(1)]
$S_p$ is cotype $(2,H)$ if and only if the formal identity $$id : RC[p] \rightarrow H$$ is completely bounded.
\item[(2)]
$L_p(\mu)$ is cotype $(2,H)$ if and only if the formal identity $$id : RC[p'] \rightarrow H$$ is completely bounded.
\end{itemize}

\end{thm}
\begin{proof}
(1) Suppose $S_p$ is cotype $(2,H)$. Then since $R[p], C[p] \subseteq S_p$ the formal identities
$$id : C[p] \rightarrow H, \;\text{and} \; id : R[p] \rightarrow H$$ are completely bounded by Proposition \ref{prop-cotype(2,H)-Hilbertian},
so that we get the desire conclusion. For the converse direction it suffices to show that $S_p$ is cotype $(2,RC[p])$, 
which is obtained from the Banach space case and the following fact. (Proposition 4.2.6. in \cite{J}) $$B(RC[p], S_p) = CB(RC[p], S_p).$$

(2) By a usual localization argument we can assume that $L_p(\mu) = L_p[0,1]$. Suppose $L_p[0,1]$ is cotype $(2,H)$. 
Note that $Rad_p \subseteq L_p[0,1]$ and $Rad_p \cong RC[p']$ completely isomorphically, where $Rad_p$ refers to the closed linear span of
the classical Rademacher system $\{r_i\}$ in $L_p[0,1]$. Thus, the formal identity $$id : RC[p] \rightarrow H$$ is completely bounded 
by Proposition \ref{prop-cotype(2,H)-Hilbertian}, so that we get the desired conclusion. For the converse direction it suffices to show that 
$L_p[0,1]$ is cotype $(2,RC[p'])$, which is obtained from the Banach space case and the following fact obtained similarly as in 
Proposition 4.2.6. in \cite{J}. $$B(RC[p'], L_p[0,1]) = CB(RC[p], L_p[0,1]).$$

\end{proof}

For a certain choice of $H$ we can recover the same behavior of type and cotype as in the Banach space case.
\begin{cor}\label{cor-Recovering-Lp}
$S_p$ ($1\leq p \leq 2$) has type $(p, R + C)$ and cotype $(2, R+C)$, 
and $S_q$ ($2\leq q < \infty$) has type $(2, R \cap C)$ and cotype $(q, R\cap C)$.
\end{cor}
\begin{proof}
First we consider the case $1\leq p\leq 2$. Since $S_2$ has type $(2,OH)$ and the formal identity 
$$id : OH \rightarrow R+C$$ is a complete contraction $S_2$ has type $(2,R + C)$. Thus, $S_p$ has type $(p, R + C)$ by interpolation and 
cotype $(2,R+C)$ by Theorem \ref{thm-cotype(2,H)-Lp}.

The case $2\leq q < \infty$ is obtained by duality.
\end{proof}

\section{Applications}

\subsection{Completely $(q,2)$-summing maps and $S_q$-cotype}\label{sec-littleGT}

Now we present our operator space version of ``generalized little Grothendieck's theorem".

\begin{thm}\label{thm2}
Let $F$ be operator spaces with $S_q$-cotype ($2\leq q <\infty)$. Then we have
$$B(C(K),F) \subseteq \Pi^o_{q,2}(C(K),F).$$
\end{thm}
\begin{proof}
Let $u \in B(C(K),F)$. Then, since $F$ is cotype $q$, we have that $u \in \Pi_{r}(E,F)$ for all $q<r<\infty$
from the Banach space result. (Theorem 11.14 in \cite{DJT})
Thus, we have by a similar calculation as in the proof of Theorem 11.13 in \cite{DJT} that
\begin{align*}
\begin{split}
\norm{(ux_{ij})}_{S_q(F)} & \leq C^o_q(F) \Big[ \int_{\Om} \norm{u(\sum_{i,j}g_{ij}(\om)x_{ij})}^2_F dP(\om) \Big]^{\frac{1}{2}}\\
& \lesssim C^o_q(F) \pi_r(u) \norm{(x_{ij})}_{S_2 \otimes_\lambda C(K)}\\& = C^o_q(F) \pi_r(u) \norm{(x_{ij})}_{S_2 \otimes_{\min}C(K)}.
\end{split}
\end{align*}

\end{proof}

\begin{rem}
{\rm 
\begin{itemize}
\item[(1)]
$S_q$-cotype conditions in Theorem \ref{thm2} are essential. Indeed, for $n\geq 1$, $2\leq q < p < \infty$
and the formal identity $I_n : \ell^n_{\infty} \rightarrow R_n[p]$ we have $\norm{I_n} = \sqrt{n}$.
As in the proof of Theorem \ref{compute1}, we get a subspace $$F (\cong R_n[p] \otimes_h R_n[q])\,\, \text{of} \,\, S^n_q(R_n[p])$$ 
isometric to $S^n_r$ under the mapping $e_{1j}\otimes e_{1i} \mapsto e_{ij}$, where $r = \frac{2pq}{pq + p -q} < 2$.
Now we have $$\norm{\sum^n_{j=1}e_j\otimes e_{1j}}_{S^n_2\otimes_{\min}\ell^n_{\infty}} =
\norm{\sum^n_{j=1}e_j\otimes e_{1j}}_{\ell^n_{\infty}(S^n_2)} = 1$$ and 
$$\norm{\sum^n_{j=1}e_{1j}\otimes e_{1j}}_{S^n_q(R_n[p])} = \norm{\sum^n_{j=1}e_{jj}}_{S^n_r} = n^{\frac{1}{r}}.$$
Consequently, we get $$\frac{\pi^o_{q,2}(I_n)}{\norm{u_n}} \geq n^{\frac{1}{r} - \frac{1}{2}}.$$

\item[(2)]
Unlike completely $p$-summing property, completely $(q,2)$-summing property ($2 < q$)
does not imply complete boundedness in general. Indeed, for $q>2$ and the formal identity $I_n : \min \ell^n_2 \rightarrow OH_n$ 
we have $$\pi^o_{q,2}(I_n) = \pi^o_{q,2}(I_{OH_n}) \leq n^{\frac{1}{4} + \frac{1}{2q}}$$ 
by Lemma 2.7 of \cite{L} and $$n \lesssim \norm{\min \ell^n_2 \stackrel{I_n}{\rightarrow} OH_n \stackrel{I_n^*}{\rightarrow} 
\max\ell^n_2}_{cb} \leq \norm{I_n}^2_{cb}$$ by Theorem 3.8 of \cite{P4}. Thus, we have 
$$\frac{\norm{I_n}_{cb}}{\pi^o_{q,2}(I_n)} \gtrsim n^{\frac{1}{4} - \frac{1}{2q}}.$$
\end{itemize}
}
\end{rem}

\subsection{An operator space version of Maurey's extension theorem}\label{sec-Maurey}

In this section we consider an operator space version of Maurey's extension theorem and Kwapie\'{n}'s theorem.
We fix a perfectly Hilbertian operator space $H$.

\begin{thm} Let $E$ and $F$ be operator spaces with type $(2,H)$ and cotype $(2,H^*)$, respectively.
Then for any subspace $G \subseteq E$ and any bounded linear map
$u : G \rightarrow F$ we have an extension
$$\tilde{u} : E \rightarrow F \,\, \text{with}\,\,\, \gamma_H(\tilde{u}) \leq T_{2,H}(E)C_{2,H^*}(F)\norm{u}.$$
\end{thm}
\begin{proof}
First we observe that we can reduce our theorem to the case that $G$ and $F$ are finite dimensional by a standard argument.
We fix $u : G \rightarrow F$ and assume that for some constant $C>0$ there are extensions
$$u_Z : E \stackrel{A_Z}{\rightarrow} H(I_Z) \stackrel{B_Z}{\rightarrow} F \,\, \text{with}\,\,\, \norm{A_Z}_{cb}\leq 1
\,\,\text{and}\,\, \norm{B_Z}_{cb}\leq C\norm{u}$$ for all finite dimensional $Z \subseteq G$.
Now we consider a nontrivial ultrafilter $\U$ of the set of all finite dimensional subspaces of $G$ ordered by inclusion.
If we set $$A : E \rightarrow \prod_{\U} H(I_Z), \,\, x\mapsto (A_Z x), \,\, B : \overline{A(G)} \rightarrow F, \,\, (A_Zx) \mapsto ux$$
then we have $$\norm{A}_{cb} \leq \lim_{\U}\norm{A_Z}_{cb} \leq 1 \,\,\text{and}\,\, \norm{B}_{cb} \leq \lim_{\U}\norm{B_Z}_{cb} \leq C\norm{u},$$
which leads us to our desired extension $\tilde{u} = BPA$, where $P$ is the orthogonal projection from $\prod_{\U}H(I_Z)$ onto $\overline{A(G)}$
since the class of $H(I)$'s for some index set $I$ is closed under ultraproduct.

Now we can assume that $G$ is finite dimensional, then since the range of $u$ is finite dimensional we can assume that so is $F$.
Let's fix $u : G \rightarrow F$, and consider any $v : F \rightarrow G$. Note that the subquadratic conditions for $H$ and $H^*$ 
together with the Remark in p.82 of \cite{P2} enable us to use Theorem 6.1 of \cite{P2} in our situation. 
Thus, by Theorem 6.1 of \cite{P2} there is a factorization
$$iv : F \stackrel{A}{\rightarrow} \ell_2(I) \stackrel{B}{\rightarrow} E\,\, \text{with} \,\, \pi_{2,H^*}(A) \leq 1\,\,
\text{and}\,\, \pi_{2, H}(B^*) \leq \gamma^*_H(iv),$$ where $i : G \hookrightarrow E$ is the inclusion.
If we set $\widetilde{B} = BPj$, where $P$ is the orthogonal projection from $H(I)$ onto $\overline{\text{ran}(A)}$ and
$j : \overline{\text{ran}(A)} \hookrightarrow H(I)$ is the inclusion, then we have the factorization
$$v : F \stackrel{A}{\rightarrow} \overline{\text{ran}(A)} \stackrel{\widetilde{B}}{\rightarrow} G,$$
so that we have
\begin{align*}
\begin{split}
\abs{\text{tr}(vu)} & = \abs{\text{tr}(\widetilde{B}Au)} \leq \ell^*(Au) \ell(\widetilde{B}) \leq \norm{u}\ell^*(A)\ell(B)\\
& \leq \norm{u}\ell^*(A)T_{2,H}(E)\pi_{2,H}(B^*) \leq T_{2,H}(E)C_{2,H^*}(F)\norm{u}\pi_{2,H^*}(A)\pi_{2,H}(B^*)\\
& \leq T_{2,H}(E)C_{2,H^*}(F)\norm{u}\gamma^*_H(iv).
\end{split}
\end{align*}
By applying the Hahn-Banach theorem to the functional $v \mapsto \text{tr}(vu)$ we can obtain the desired extension
$$\tilde{u} : E \rightarrow F \,\, \text{with}\,\,\, \gamma_H(\tilde{u}) \leq T_{2,H}(E)C_{2,H^*}(F)\norm{u}.$$

\end{proof}

\begin{cor} Every operator space with type $(2,H)$ and cotype $(2,H^*)$ is completely
isomorphic to $H(I)$ for some index set $I$.
\end{cor}

\bibliographystyle{amsplain}
\providecommand{\bysame}{\leavevmode\hbox
to3em{\hrulefill}\thinspace}

\end{document}